\documentclass{amsart}
\usepackage{setspace}

\newtheorem{theorem}{Theorem}[section]
\newtheorem{lemma}[theorem]{Lemma}
\newtheorem{proposition}[theorem]{Proposition}

\theoremstyle{definition}

\theoremstyle{remark}

\numberwithin{equation}{section}

\begin{document}

\title{Transitivity and Bundle Shifts}

\author{Ronald G. Douglas}
\address{Department of Mathematics, Texas A\&M University, College Station, Texas 77843-3368}
\email{rgd@tamu.edu}

\author{AnJian Xu}
\address{School of Mathematical Sciences and Statistics,
Chongqing University of Technology, Chongqing, China 400075}
\email{xuaj@cqut.edu.cn}
\curraddr{Department of Mathematics, Texas A\&M University, College Station, Texas 77843-3368}
\thanks{The second author was supported in part by NSFC Grant \#11271388.}


\date{February 28, 2014.}

\dedicatory{ }


\begin{abstract}
A subalgebra $A$ of the algebra $B(\mathcal{H})$ of bounded linear operators on a separable Hilbert space $\mathcal{H}$ is said to be catalytic if every transitive subalgebra $\mathcal{T}\subset B(\mathcal{H})$ containing it is strongly dense. We show that for a hypo-Dirichlet or logmodular algebra, $A=H^{\infty}(m)$ acting on a generalized Hardy space $H^{2}(m)$ for a representing measure $m$ that defines a reproducing kernel Hilbert space is catalytic.\\
\indent
For the case of a nice finitely-connected domain, we show that the "holomorphic functions" of a bundle shift yields a catalytic algebra, thus generalize a result of Bercovici, Foias, Pearcy and the first author[7].
\end{abstract}

\maketitle

\section{Introduction}
\indent
In this  paper, let $\mathcal{H}$ denote a complex separable Hilbert space, and $B(\mathcal{H})$ the algebra of all bounded linear operators on $\mathcal{H}$. A unital subalgebra $\mathcal{A}$ of $B(\mathcal{H})$ is called \emph{transitive} if it has only trivial invariant subspaces. The transitive algebra problem asks if every transitive algebra $\mathcal{A}\subset B(\mathcal{H})$ is strongly dense in $B(\mathcal{H})$. An operator, or a set of operators, will be called \emph{catalytic} if any transitive operator \linebreak algebra containing it is strong dense in $B(\mathcal{H})$. It is Arveson [6] who first stated explicitly the problem, and developed the main tool for \linebreak studying transitive algebras. In the same paper, he proved that both the unilateral shift with multiplicity one and a non scalar normal operator of multiplicity one are catalytic. Later, Richter [14] proved that the Dirichlet shift is catalytic, and Nordgren [12,13] generalized Arveson's result to unilateral shifts with finite multiplicities. Cheng, Guo and Wang [8] proved that the coordinate multiplication operators on functional Hilbert spaces with complete Nevanlinna-Pick kernels are catalytic. The invariant subspace problem asks if a singly-generated algebra acting on a separable Hilbert space $\mathcal{H}$ can be transitive. Catalytic operators must have nontrivial invariant subspaces.\\
\indent
Let $R$ be a finitely-connected, bounded domain in the complex plane whose boundary $\partial R$ consists of $n+1$ nonintersecting smooth Jordan curves. The Hardy space $H^{2}(R)$ over $R$ is defined to be the space of analytic functions $f$ on $R$ such that the subharmonic function $|f|^{2}$ is majorized by a harmonic function $u$ on $R$. For each fixed point $t\in R$, there is a norm $\|\cdot\|_{t}$ on $H^{2}(R)$ defined by
\[\|f\|_{t}=\inf\limits_{u}\{u(t)^{\frac{1}{2}}|\hspace{2pt}u\text{ is a harmonic majorant of }|f|^{2}\}.\]
\indent
As in the case of $R$ equal to the unit disc $\mathbb{D}$, one can define an isomorphic spaces of functions on $\partial R$ using boundary values. Let $m_{t}$ be the harmonic measure for the point $t\in R$; that is, $f(t)=\int_{\partial R}fdm_{t}$ for $f$ bounded and harmonic on $R$. Let $L^{2}(\partial R,m_{t})$ be the space of square-integrable complex-valued measurable functions on $\partial R$ defined with respect to $m_{t}$. Then $H^{2}(\partial R,m_{t})$ is defined to be the set of functions $f\in L^{2}(\partial R,m_{t})$ such that $\int_{\partial R}f(z)g(z)dm_{t}=0$ for every $g$ that is analytic in a neighborhood of the closure of $R$. Although the norm on $H^{2}(\partial R,m_{t})$ depends on the fixed point $t\in R$, the spaces of functions for both $H^{2}(\partial R,m_{t})$ and $L^{2}(\partial R,m_{t})$ are independent of $t$, and these norms are boundedly equivalent for all $t\in R$. We fix $t\in R$ and use $H^{2}(\partial R)=H^{2}(\partial R,m_{t})$ omitting reference to it. \\
\indent
For any point $w\in R$, the point evaluation functional on $H^{2}(R)$, defined by $ev_{w}(f)=f(w)$, is a bounded linear functional. Thus $H^{2}(R)$ is a reproducing kernel Hilbert space for $R$. We use $k_{w}$ to denote the reproducing kernel function for $w\in R$; that is, $ev_{w}(f)=\langle f, k_{w}\rangle$. Using the isomorphism of $H^{2}(R)$ and $H^{2}(\partial R)$, we see that $H^{2}(\partial R)$ is a reproducing kernel Hilbert space also. Let $\hat{k}_{\lambda}\in H^{2}(\partial R)$ denote the reproducing kernel function for $H^{2}(\partial R)$ at $\lambda\in R$; that is, $\langle f, \hat{k}_{\lambda}\rangle=f(\lambda)$ for $f\in H^{2}(\partial R)$. (Of course, for the identification of $H^{2}(R)$ and $H^{2}(\partial R)$ to be an isomorphism, one must use the same point $t\in R$ in defining both norms. One can also use the function $w\rightarrow\langle g,\hat{k}_{w}\rangle$ for $g\in H^{2}(\partial R)$ to identify $H^{2}(R)$ and $H^{2}(\partial R)$.)\\
\indent
We define an operator $T_{z}$ on $H^{2}(R)$ by $(T_{z}f)(z)=zf(z)$ for every $f\in H^{2}(R)$ and $z\in R$. Furthermore, the isomorphism of $H^{2}(R)$ and $H^{2}(\partial R)$ induces a corresponding operator, also denoted $T_{z}$, on $H^{2}(\partial R)$. And we define $N_{z}$ on $L^{2}(\partial R)$ by  $(N_{z}f)(z)=zf(z)$ for $z\in\partial R$. It can be seen that $T_{z}$ is a pure subnormal operator [3] and $N_{z}$ is the minimal normal extension of $T_{z}$ on $H^{2}(\partial R)$. Let $T_{\varphi}$ be the operator on $H^{2}(R)$ defined by $(T_{\varphi}f)(z)=\varphi(z)f(z)$ for every $f\in H^{2}(R)$ and $\varphi\in H^{\infty}(R)$. We use $T_{\varphi}$ to denote the corresponding operator on $H^{2}(\partial R)$.\\
\indent
Similarly, for a Hilbert space $\mathcal{H}$, we can define an $\mathcal{H}$-valued Hardy \linebreak space, $H_{\mathcal{H}}^{2}(R)$, which is the space of $\mathcal{H}$-valued analytic functions \linebreak $f:R\rightarrow\mathcal{H}$ such that the subharmonic function $\|f(z)\|^{2}_{\mathcal{H}}$ is majorized by a harmonic function $u$ on $R$. We define two corresponding operators, $(T_{\mathcal{H}}f)(z)=zf(z)$ for $f\in H_{\mathcal{H}}^{2}(R)$ and $z\in R$, and $N_{\mathcal{H}}$ on $L_{\mathcal{H}}(\partial R)$, $(N_{\mathcal{H}}f)(z)=zf(z)$ for $f\in L_{\mathcal{H}}^{2}(\partial R)$ and $z\in\partial R$. Now $H_{\mathcal{H}}^{2}(R)$ is a reproducing kernel Hilbert space, and we use $k_{\lambda}^{\mathcal{H}}\in B(\mathcal{H})$ to represent the reproducing kernel at $\lambda\in R$; that is, $\langle f(\lambda),h\rangle_{\mathcal{H}}=\langle f,k_{\lambda}^{\mathcal{H}}h\rangle_{H_{\mathcal{H}}^{2}(R)}$ for $f\in H_{\mathcal{H}}^{2}(R)$ and $h\in\mathcal{H}$. For more information about function theory on finitely connected domains, one can see [1,15,16]. \\

\indent Let $E$ be a Hermitian holomorphic vector bundle over $R$. A section of $E$ is a holomorphic function $f$ from $R$ into $E$ such that $p(f(z))=z$ for all $z\in R$, where $p:E\rightarrow R$ is the projection map [17]. The set of all holomorphic sections of $E$ is denoted by $\Gamma_a(E)$. For $E$ a rank $n$ Hermitian vector bundle over $R$, a unitary coordinate cover for $E$ is a covering $\{U_{s},\varphi_{s}\}$ with $\varphi_{s}:U_{s}\times\mathbb{C}^{n}\rightarrow E|_{U_{s}}$ such that for each $s$ and $z\in U_{s}$, the fiber map $\varphi_{s}^{z}:\mathbb{C}^{n}\rightarrow E_{z}$, is unitary. The unitary coordinate cover $\{U_{s},\varphi_{s}\}$ is said to be \emph{flat} if the transition functions, $\varphi_{st}=\varphi_{s}^{-1}\varphi_{t}$ on $U_{s}\cap U_{t}$ for all $s$ and $t$, are constant. A flat unitary vector bundle is a vector bundle with a flat unitary coordinate covering.\\

\indent
If $E$ is a flat unitary vector bundle over the finitely-connected domain $R$ with fiber $\mathcal{E}$ and coordinate covering $\{U_{s},\varphi_{s}\}$ and $f$ is a holomorphic section of $E$, then for $z\in U_{s}\cap U_{t}$, the operator $(\varphi_{t}^{z})^{-1}\varphi_{s}^{z}$ is unitary so that $\|(\varphi_{t}^{z})^{-1}f(z)\|=\|(\varphi_{s}^{z})^{-1}f(z)\|$. This means that there is a function on $R$ defined by $h_{f}^{E}(z)=\|(\varphi_{s}^{z})^{-1}f(z)\|_{E}$, where $z\in U_{s}$. One defines $H_{E}^{2}(R)$ to be the space of holomorphic sections $f$ of $E$ such that $(h_{f}^{E}(z))^{2}$ is majorized by a harmonic function, then $H_{E}^{2}(R)$ is a Hilbert space. $H_{E}^{2}(R)$ is invariant under multiplication by any bounded analytic function on $R$. The operator $T_{E}$ on $H_{E}^{2}(R)$, defined by $(T_{E}f)(z)=zf(z)$ for $z\in R$, is called a \emph{bundle shift} over $R$. These objects are studied by Abrahamse and the first author [3].\\

\indent
Finally, let $\mathcal{T}_{E}$ be the subalgebra of $B(H_{E}^{2}(R))$ of operators $T_{\Phi}$, where $\Phi$ is a bundle map on $E$ which extends to an open set containing the closure of $R$.\\

\indent
In this paper, first for a hypo-Dirichlet algebra, we prove that an operator algebra on $H^{2}(m)$ containing $H^{\infty}(m)$ is strongly dense in $B(H^{2}(m))$, which implies that any operator algebra containing $\mathcal{Q}(T_{z})$ is strongly dense in $B(H^{2}(R))$, where $\mathcal{Q}$ is the algebra consisting of all rational analytic functions with poles outside the closure of $R$, Thus this algebra is catalytic. We next extend the result to the context of logmodular algebras. Second, we show that an analogous catalytic result also holds for the case of bundle shifts, which generalizes the result of Bercovici, Foias, Pearcy and the first author [7]. We do this by proving that the algebra $\mathcal{A}\otimes M_{n}(\mathbb{C})$ is catalytic on $\mathcal{K}\otimes\mathbb{C}^{n}$ if $\mathcal{A}$ is catalytic on $\mathcal{K}$ for any positive integer $n$. This fact was observed by Nordgren [13].\\
\indent
We conclude this introduction by recalling that these techniques don't work for the Bergman shift, and hence are limited in scope [9].\\
\indent
\textbf{Acknowledgments}. Research of the second author was carried out during a year-long visit to the Department of Mathematics at Texas A\&M University, which was supported by the China scholarship council and the Department of Mathematics at Chongqing University of Technology.

\section{Function Algebras}
\indent
Let $X$ be a compact Hausdorff space and $C(X)$ be the algebra of all continuous complex-valued functions on $X$. A subalgebra $A$ of $C(X)$ is called a \emph{uniform algebra} on $X$ if it contains the constant functions and separates the points of $X$. One can see Gamelin's book [10] for details on uniform algebras. The space of real parts of the functions in $A$ is denoted by $ReA$, the set of invertible elements in $A$ by $A^{-1}$, and $\log|A^{-1}|$ denotes the set, $\{\log|f|\hspace{2pt}, f\in A^{-1}\}$.\\
\indent
Let $\varphi$ be a multiplicative linear functional on $A$. Then there exist positive measures $m$ on $X$ such that
\[\varphi(f)=\int_{X}fdm,\hspace{6pt}\text{for }f\in A,\]
which are called \emph{representing measures} for $\varphi$. Usually, the representing measure is not unique. An \emph{Arens-Singer measure} for $\varphi$ is a positive measure $m$ on $X$ such that
\[\log|\varphi(f)|=\int_{X}\log|f|dm,\hspace{6pt}\text{for }f\in A^{-1}.\]
An Arens-Singer measure for $\varphi$ is also a representing measure for $\varphi$. In fact, we have
\begin{eqnarray}
&\int_{X}fdm&=\int_{X}(u+iv)dm=\int_{X}\log e^{u}dm+i\int_{X}\log e^{v}dm\nonumber\\
&&=\log e^{\varphi(u)}+i\log e^{\varphi(v)}=\varphi(u)+i\varphi(v)=\varphi(f)\nonumber
\end{eqnarray}
for $f=u+iv$. An Arens-Singer measure always exists (cf. [4,11] and the references therein).

\subsection{Hypo-Dirichlet Algebras}
The uniform algebra $A$ is called a \emph{hypo-Dirichlet algebra} if there exists a finite set of elements $\{f_{1},\cdots,f_{s}\}$ in $A^{-1}$ such that the linear span of $ReA$ and $\log|f_{1}|,\cdots,\log|f_{s}|$ is \linebreak uniformly dense in $ReC(X)$, where s is the minimal number possible. Ahern and Sarason studied hypo-Dirichlet algebras in [4].\\
\indent
For a multiplicative functional $\varphi$ on $A$, we use $\mathcal{M}(\varphi)$ to denote the set of representing measures for $\varphi$, and $S$ to denote the real linear span of all differences between pairs of measures in $\mathcal{M}(\varphi)$. For a hypo-Dirichlet algebra, one can show that $\dim S=s$.
Further, there is only one Aren-Singer measure on a hypo-Dirichlet algebra.\\

\indent
For a fixed measure $m$ on $X$ corresponding to the multiplicative linear functional $\varphi$ on $A$, let $A_{0}$ denote the kernel of the functional $\varphi$, $H^{2}(m)$ denote the closure of $A$ in $L^{2}(m)$, and let $H^{\infty}(m)$ be the weak-star closure of $A$ in $L^{\infty}(m)$. Ahern and Sarason [4] showed that the orthogonal complement $N$ of $A+\bar{A}$ in $L^{2}(m)$ is a subspace of $L^{\infty}$. Moreover $L^{2}(m)$ has the following decomposition, $L^{2}(m)=H^{2}(m)\oplus \overline{H^{2}_{0}(m)}\oplus N$, where $\overline{H^{2}_{0}(m)}$ is the complex conjugates of the functions in $H^{2}_{0}(m)$, the functions which are annulled by $m$.\\

\indent
If the difference of two multiplicative functionals $\varphi$, $\psi$ on $A$ has \linebreak norm less than 2, then we say that $\varphi$ and $\psi$ lie in the same \emph{part} of the maximal ideal space $\mathcal{M}_{A}$ of $A$, which defines an equivalence relation on $\mathcal{M}_{A}$ and partitions $\mathcal{M}_{A}$ into the disjoint union of parts. In the case of a hypo-Dirichlet algebra, the Arens-Singer measures of two functionals in the same part are mutually, boundedly absolutely continuous [10].\\

\indent
An element $x$ of a convex subset $K$ of a vector space $V$ is said to be a \emph{core point} of $K$ if whenever $y\in V$ is such that $x+y\in K$, then $x-\epsilon y\in K$ for $\epsilon>0$ sufficiently small. A measure $m$ is called a \emph{core measure} for $\varphi$ if $m$ is a core point of $\mathcal{M}(\varphi)$. In particular, one can show that the harmonic measure on $\partial R$, corresponding to each $z\in R$, is a core measure for $\varphi$.\\

\indent
With these notions defined, we can state a basic lemma on representations of functions in a generalized Hardy space by combining Theorem V.4.2 and Theorem V.4.3 in Gamelin[10]. The interested reader can consult Gamelin [10] for further information.

\begin{lemma}
Let $A$ be a function algebra on a compact metric space $X$, and $\varphi$ be a multiplicative linear functional on $A$. If the space of representing measures for $\varphi$ is finite dimensional and if $m$ is a core measure, then every function $f\in H^{2}(m)$ can be expressed as a quotient of two $H^{\infty}(m)$ functions
\end{lemma}

\subsection{Logmodular Algebras}
\indent
A uniform algebra $A$ is a \emph{Dirichlet algebra} on $X$ if $ReA$ is uniformly dense in $ReC(X)$; $A$ is a \emph{logmodular algebra} on $X$ if $\log|A^{-1}|$ is uniformly dense in $ReC(X)$. A Dirichlet algebra is a logmodular algebra since $ReA\subseteq\log|A^{-1}|$. One can consult Hoffman[10] for further results about logmodular algebras.\\
\indent
For every multiplicative linear functional $\varphi$ on a logmodular algebra $A$, there is a unique representing measure $m$ on $X$. Let $A_{0}$ denote the kernel of the functional $\varphi$ also, $H^{2}(m)$ denote the closure of $A$ in $L^{2}(m)$, and $H^{\infty}(m)$ be the weak-star closure of $A$ in $L^{\infty}(m)$. A function $g$ in $H^{2}(m)$ is called an \emph{outer function} if $Ag$ is dense in $H^{2}(m)$. Equivalently, $g$ is outer if $\int\log|g|dm=\log|\int gdm|>-\infty$.\\
\indent
Now we have the following lemma whose proof is similar to that of the classic Hardy space $H^{2}(\mathbb{D})$ on $\mathbb{D}$ [12], and Theorem V.4.3. in Gamelin[10].
\begin{lemma}
Let $A$ be a logmodular algebra on $X$. For every multiplicative linear functional $\varphi$ on $A$, let $m$ be the unique representing measure for $\varphi$ on $X$. For every function $f\in H^{2}(m)$, there exist two functions $f_{1}\in H^{\infty}(m)$, $f_{2}\in H^{\infty}(m)^{-1}$ such that $f=\frac{f_{1}}{f_{2}}$.
\end{lemma}
\begin{proof}
We have the orthogonal decomposition [11] of the space $L^{2}(m)=H^{2}(m)\oplus\overline{H_{0}^{2}(m)}$, where $H_{0}^{2}(m)$ is the closure of $A_{\varphi}$ in $H^{2}(m)$. Note that for every function $u\in Re H^{2}(m)$, there exists a unique function $v\in H^{2}_{0}(m)$ such that $u+iv\in H^{2}(m)$. Therefore, for $f\in H^{2}(m)$, there  exists a function $v$ such that $\log^{+}|f|+iv\in H^{2}(m)$ where $\log^{+}|f|=\max\{\log|f|,0\}$. If we define $g=e^{-(\log^{+}|f|+iv)}$, then $g$ is an outer function and invertible in $H^{\infty}(m)$, hence $fg$ is bounded also. So we have the representation lemma for functions in $H^{2}(m)$ by setting $f_{1}=fg$ and $f_{2}=g$.
\end{proof}
\indent
This Lemma generalizes the corresponding result for the classic Hardy space $H^{2}(\mathbb{D})$ on the unit disc.

\section{Catalytic in Function algebras }
\indent
Let $\mathcal{H}$ denote an infinite dimensional Hilbert space, and $B(\mathcal{H})$ be the algebra of bounded linear operators on $\mathcal{H}$. Let $\mathcal{A}$ be a subalgebra of $B(\mathcal{H})$, and let $n$ be a positive integer. Then $\mathcal{A}$ is called \emph{$n$-transitive} if for any choice of elements $x_{1},\cdots,x_{n},y_{1},\cdots,y_{n}\in \mathcal{H}$ with $\{x_{1},\cdots,x_{n}\}$ linearly independent, there exists a sequence $A_{i}\in\mathcal{A}$ such that
\[\lim\limits_{j}A_{j}x_{k}=y_{k},\hspace{6pt}1\leq k\leq n.\]
Further, a unital subalgebra $\mathcal{A}$ of $B(\mathcal{H})$ is said to be \emph{transitive} if it is 1-transitive or, equivalently, it has only trivial invariant subspaces; that is, only $\{0\}$ and $\mathcal{H}$. It is easy to show that a subalgebra $\mathcal{A}$ is strongly dense in $B(\mathcal{H})$ if and only if $\mathcal{A}$ is $n$-transitive for every $n\geq1$. This equivalence is due essentially to Arveson [6].\\
\indent
As mentioned in the introduction, Arveson introduced an inductive strategy to show that transitive algebras are strongly dense. The first step in this approach is the content of the following lemma.
\begin{lemma}
(Arveson's Corollary 2.5.[6]) Let $\mathcal{A}$ be a transitive subalgebra of $B(\mathcal{H})$, where $\mathcal{H}$ is a Hilbert space. Then $\mathcal{A}$ is not 2-transitive, if and only if, there exists a closed, densely defined, non-scalar linear transformation on $\mathcal{H}$ that commutes with $\mathcal{A}$.
\end{lemma}
\indent
The further induction steps are all the same and are subsumed in the following lemma.
\begin{lemma}
(Arveson's Corollary 2.5.[6]) Let $A$ be a subalgebra of $B(\mathcal{H})$ for a Hilbert space $\mathcal{H}$. For $N\geq2$, suppose that $\mathcal{A}$ is $N$-transitive but not $(N+1)$-transitive. Then there exists $N$ linear transformations $T_{1},\cdots,T_{N}$, defined on a common linear domain $\mathcal{D}$ dense in $\mathcal{H}$, such that\\
\indent\hspace{6pt}(i) each $T_{i}$ commutes with $\mathcal{A}$;\\
\indent\hspace{6pt}(ii) no $T_{i}$ is closable; (that is, the closure of the graph of $T_{i}$ is not the graph of a transformation); and\\
\indent\hspace{6pt}(iii) $\{(x,T_{1}x,\cdots,T_{N}x)|\hspace{2pt}x\in\mathcal{D}\}$ is a proper closed subspace of $\mathcal{H}^{(N+1)}$.
\end{lemma}
\indent
For a general hypo-Dirichlet algebra $A$, a part of $M_{A}$ is defined as a analytic part if it contains more than one point. Firstly we use Lemma 2.1 to  prove the following proposition.
\begin{proposition}
Let $A$ be a general hypo-Dirichlet algebra on $X$, $\varphi$ be a multiplicative linear functional on $A$, and assume that the representing measure $m$ corresponding to $\varphi$ is a core measure, and $H^{2}(m)$ is a reproducing kernel Hilbert space. If $\mathcal{A}$ is an operator algebra containing
\[\{T_{h}\in B(H^{2}(m))|\hspace{4pt}h\in H^{\infty}(m)\},\]
$T$ is a map from a dense $\mathcal{A}$-invariant submanifold $\mathcal{D}$ of $H^{2}(m)$ to $H^{2}(m)$, and suppose $TAf=ATf$,$f\in\mathcal{D}$, $A\in \mathcal{A}$. Then there exists a function of the form $h=\frac{f_{1}}{f_{2}}$, with $f_{1},f_{2}\in H^{\infty}(m)$, such that $Tf=hf$ for all $f\in\mathcal{D}$. Moreover, $T$ is linear and closable.
\end{proposition}
\begin{proof}
\indent
Firstly, we prove that $\frac{Tf}{f}$ is equal to $\frac{Tg}{g}$ for every pair $f,g\in\mathcal{D}$, $f,g\neq0$, the zero function. By the statement of Lemma 2.1, there exist functions $f_{1},f_{2},g_{1},g_{2}\in H^{\infty}(m)$, $f_{2}$ and $g_{2}$ invertible in $H^{\infty}(m)$ such that $f=\frac{f_{1}}{f_{2}}$ and $g=\frac{g_{1}}{g_{2}}$. We want to prove $\frac{Tf}{f}=\frac{Tg}{g}$, or
$g_{1}f_{2}Tf=g_{2}f_{1}Tg$. But from the commutativity of $T$ with $T_{k}$ for any $k\in H^{\infty}(m)$, it follows for $h$, the quotient $\frac{Tf}{f}$, for any $f\in\mathcal{D}$ such that zero is not in the range of $f$, we have $Tf=hf$. If $h$ is the quotient of two functions in $H^{2}(m)$, then since every function in $H^{2}(m)$ can be expressed as a quotient of two functions in $H^{\infty}(m)$. Hence $h$ can be written as a quotient $\frac{h_{1}}{h_{2}}$ of two functions $h_{1},h_{2}\in H^{\infty}(m)$. \\
\indent
Secondly, we prove that the operator defined by $Tf=hf$ on $\mathcal{D}=\{f\in H^{2}(m)|\quad hf\in H^{2}(m)\}$ is closable. It is clear that the domain $\mathcal{D}$ is invariant under linear operations in $H^{2}(m)$. Further for any sequence $\{f_{n}\}\subset\mathcal{D}$ such that $\lim\limits_{n\rightarrow\infty}f_{n}=f$ and $\lim\limits_{n\rightarrow\infty}hf_{n}=g$, $f,g\in H^{2}(m)$, we see that $g=hf$ as follows. We assume that the reproducing kernel for $\lambda$ in analytic parts of $M_{A}$ is denoted by $k^{m}_{\lambda}$, $\lim\limits_{n\rightarrow\infty}f_{n}(\bar{\lambda})=\lim\limits_{n\rightarrow\infty}\langle f_{n},k^{m}_{\lambda}\rangle=\langle f,k^{m}_{\lambda}\rangle=f(\bar{\lambda})$, and $\lim\limits_{n\rightarrow\infty}f_{n}h(\bar{\lambda})=\lim\limits_{n\rightarrow\infty}\langle hf_{n},k^{m}_{\lambda}\rangle=\langle g,k^{m}_{\lambda}\rangle=g(\bar{\lambda})$. So $g(\bar{\lambda})=h(\bar{\lambda})f(\bar{\lambda})$. This means that $hf\in H^{2}(m)$ and $g=hf$; that is, $T$ is closable.
\end{proof}
\begin{lemma}
Let T be the operator defined in Proposition 3.3, $k^{m}_{\lambda}(z)$ be the kernel function at the point $\lambda$ in analytic parts of $M_{A}$. Then $k^{m}_{\lambda}$ is in the domain of the operator $T^{\ast}$.
\end{lemma}
\begin{proof}
For any $f\in\mathcal{D}$,
\[\langle Tf,k_{\lambda}\rangle\langle hf,k^{m}_{\lambda}\rangle=h(\bar{\lambda})f(\bar{\lambda})=h(\bar{\lambda})\langle f,k^{m}_{\lambda}\rangle=\langle f,\overline{h(\bar{\lambda})}\rangle=\langle f,T^{\ast}k^{m}_{\lambda}\rangle,\]
it is continuous in $f$. This proves that $k^{m}_{\lambda}$ is in the domain of $T^{\ast}$ and $T^{\ast}k^{m}_{\lambda}=\overline{h(\lambda)}k^{m}_{\lambda}$.
\end{proof}
\begin{theorem}
Let $A$ be a general hypo-Dirichlet algebra, and $\varphi$ be a multiplicative linear functional. If the representing measure $m$ corresponding to $\varphi$ is a core measure, and $H^{2}(m)$ is a reproducing kernel Hilbert space. Then $H^{\infty}(m)$ is catalytic on $H^{2}(m)$.
\end{theorem}
\begin{proof}
Suppose $\mathcal{A}$ is a transitive subalgebra of $B(H^{2}(m))$ containing $H^{\infty}(m)$ and $T$ is a closed, densely defined operator commuting with $\mathcal{A}$. Then $T^{\ast}$ commute with $\mathcal{A}^{\ast}$. For $\lambda$ in analytic parts of $M_{A}$, by Lemma 3.4, $k^{m}_{\lambda}$ is in the domain of $T^{\ast}$ and $T^{\ast}k^{m}_{\lambda}=\overline{h(\lambda)}k^{m}_{\lambda}$, so that for every $A\in\mathcal{A}$,
\[\overline{h(\lambda)}A^{\ast}k^{m}_{\lambda}=A^{\ast}\overline{h(\lambda)}k^{m}_{\lambda}=A^{\ast}T^{\ast}k^{m}_{\lambda}=T^{\ast}A^{\ast}k^{m}_{\lambda}.\]
This shows that the nonempty submanifold $\{f\in\mathcal{D}(T^{\ast})|\hspace{2pt}T^{\ast}f=\overline{h(\lambda)}f\}$ is invariant under $\mathcal{A}^{\ast}$. Since $\mathcal{A}^{\ast}$ is transitive, this manifold is dense. Thus $T^{\ast}=\overline{h(\lambda)}I$, and so $T$ is a scalar. This shows that  $\mathcal{A}$ is $2$-transitive.\\
\indent
Every densely defined linear transformation that commutes with $\mathcal{A}$ is closable by Proposition 3.3, hence by Arveson's Lemma 3.2, $\mathcal{A}$ is $N$-transitive for every $N>2$, which shows that $\mathcal{A}$ is strongly dense in $B(H^{2}(m))$. Therefore $H^{\infty}(m)$ is catalytic.
\end{proof}
\indent
For logmodular algebras, we can use Lemma 2.2 to prove in a similar way the following theorem by proving analogues of Proposition 3.3, Corollary 3.4 and finally Theorem 3.5 to reach the following result.
\begin{theorem}
Let $A$ be a general logmodular algebra, $\varphi$ be a multiplicative linear functional, $m$ be the representing measure corresponding to $\varphi$. If $H^{2}(m)$ is a reproducing kernel Hilbert space, then $H^{\infty}(m)$ is catalytic on $H^{2}(m)$.
\end{theorem}
\indent
Note that the domain for $H^{2}(m)$ can be identified with a subset of the maximal ideal space for $H^{\infty}(m)$.

\section{Finitely-connected domains}
\indent
We now return to the case of the finitely-connected domain $R\subset\mathbb{C}$ whose boundary $\partial R$ consists of $(n+1)$ nonintersecting smooth Jordan curves. Here we can be more precise. In particular, $H^{\infty}(R)$ can be shown to be a hypo-Dirichlet algebra on $\partial R$ with $s=n$. For any point $z\in R$, the harmonic measure $m(=m_{z})$, supported on $\partial R$, corresponding to evaluation at $z$ is the unique Arens-Singer measure for the functional $\varphi$ on $H^{\infty}(R)$ of evaluation at $z$. For any two points $z_{1},z_{2}\in R$, the corresponding functionals lie in the same part, and they are mutually boundedly absolutely continuous. Moreover, in this case, the measures in $S$ can be described explicitly. In fact, a basis for the homology of $R$ can be used to give a linear basis for $S$ [4].
\begin{proposition}
Let $R$ be a finitely-connected domain whose boundary, $\partial R$, consists of $n+1$ nonintersecting smooth Jordan curves. If $f\in H^{2}(R)$, then $f=\frac{g}{h}$, where $g$ and $h$ belong to $H^{\infty}(R)$, and $h$ is invertible in $H^{\infty}(R)$.
\end{proposition}
\begin{proof}
Let $A$ be the algebra of all function on $\partial R$ that can be uniformly approximated by rational functions with poles off $R$. It is well known that $A$ is a hypo-Dirichlet algebra [4]. For a point $z\in R$, the corresponding evaluation functional $\varphi$ is a multiplicative linear functional. The space of representing measures $M_{\varphi}$ has dimension $n$. So there exists a core measure denoted by $m$ corresponding to $\varphi$ [10, p106]. Then the proposition follows from Proposition 3.3.
\end{proof}
\indent For finitely-connected domains and the algebra $A$ in Proposition 4.1, the analytic part of $M_{A}$ is the open set $R$. So we have the following proposition which can be proved just as Proposition 3.3.
\begin{proposition}
For the finitely-connected domain $R$ whose boundary $\partial R$ consists of $(n+1)$ nonintersecting smooth Jordan curves, let $\mathcal{A}$ be an operator algebra containing
\[\{T_{h}\in B(H^{2}(R))|\hspace{4pt}h\in H^{\infty}(R)\}.\]
If $T$ is a map of a dense $\mathcal{A}$-invariant submanifold $\mathcal{D}$ of $H^{2}(R)$ to $H^{2}(R)$, such that $TAf=ATf$, for $f\in\mathcal{D}$, $A\in \mathcal{A}$. Then there exists a meromorphic function of the form $h=\frac{f_{1}}{f_{2}}$, with $f_{1},f_{2}\in H^{\infty}(R)$, such that $Tf=hf$ for all $f\in\mathcal{D}$. Moreover, $T$ is linear and closable.
\end{proposition}
\indent The following lemma follows in the same way as Lemma 3.4.
\begin{lemma}
Let T be the closed operator defined in Proposition 4.2, and $k_{\lambda}(z)$ be the kernel function at the point $\lambda\in R$. Then $k_{\lambda}$ is in the domain of the operator $T^{\ast}$.
\end{lemma}
Now we are ready to obtain the catalytic result for finitely-connected domains by combining the method of proof in Theorem 3.5 and Proposition 4.1, Proposition 4.2, and Lemma 4.3.
\begin{theorem}
Let $R$ be a finitely-connected domain whose boundary $\partial R$ consists of $n+1$ nonintersecting smooth Jordan curves. Then the algebra $\mathcal{Q}(T_{z})$ on $H^{2}(\partial R)$ is catalytic, where $\mathcal{Q}$ is the space of all rational analytic functions with poles outside the closure of $R$.
\end{theorem}
\indent
In [3] Abrahamse and the first author proved that if $E$ is a flat unitary vector bundle over a finitely-connected domain $R$, whose boundary $\partial R$ consists of $n+1$ nonintersecting smooth Jordan curves, with fiber $\mathcal{H}$, then the bundle shift $T_{E}$ is similar to $T_{\mathcal{H}}$. Thus to show that $T_{E}$ is catalytic, it is sufficient to prove the vector-valued case or that $T_{\mathcal{H}}$ is catalytic.\\
\indent
Let $\mathcal{T}_{E}(R)$ be the subalgebra of $B(H_{E}^{2}(R))$ of operators $T_{\Phi}$, where $\Phi$ is a bundle map on $E$ which extends to an open set containing the closure of $R$. Now, in particular, we will prove $\mathcal{J}_{E}(R)$ is catalytic by first showing the following general result.
\begin{proposition}
Let $\mathcal{K}$ be a Hilbert space. If $\mathcal{A}$ is a catalytic algebra on $\mathcal{K}$, then $\mathcal{A}\otimes M_{n}(\mathbb{C})$ is a catalytic algebra on $\mathcal{K}\otimes\mathbb{C}^{n}$ for any positive integer $n$.
\end{proposition}
\begin{proof}
Suppose that the transitive algebra $\mathcal{T}$ on $\mathcal{K}\otimes\mathbb{C}^{n}$ contains $\mathcal{A}\otimes M_{n}(\mathbb{C})$. We use $\{e_{k}\}$ to denote the orthonormal basis of $\mathbb{C}^{n}$. Let $E_{ij}$ be the matrix unit operator $e_{i}\otimes e_{j}$ on $\mathbb{C}^{n}$, and $\widehat{E}_{ij}=I_{\mathcal{K}}\otimes E_{ij}$.\\
\indent
Let $T$ be a densely defined linear transformation on $\mathcal{D}\subset\mathcal{K}\otimes\mathbb{C}^{n}$ which commutes with $\mathcal{T}$. Note that $\mathcal{A}\otimes M_{n}(\mathbb{C})$ contains $\{\widehat{E}_{ij}\}$, then the operator $\widehat{E}_{jj}T\widehat{E}_{ii}$ can be seen as a densely defined linear transform from $\mathcal{K}\otimes\mathbb{C}e_{i}$ to $\mathcal{K}\otimes\mathbb{C}e_{j}$ which are both isomorphic to $\mathcal{K}$. Then it is easy to see that $\widehat{E}_{jj}T\widehat{E}_{ii}$ commutes with $\mathcal{A}$. Further, $\widehat{E}_{jj}\mathcal{T}\widehat{E}_{ii}$ is transitive on $\mathcal{K}$ since $\mathcal{T}$ is transitive and contains $\mathcal{A}$. Thus $\widehat{E}_{jj}\mathcal{T}\widehat{E}_{ii}$ is strong dense in $B(\mathcal{K})$ since $\mathcal{A}$ is catalytic, which shows that $\widehat{E}_{jj}T\widehat{E}_{ii}$ is a constant $a_{ij}\in\mathbb{C}$. So $T$ has a form $I_{\mathcal{K}}\otimes A$, where $A=(a_{ij})$ is an operator on $\mathbb{C}^{n}$, it shows that $T$ is closable. By Arveson's Lemma 3.2, $\mathcal{T}$ is $n$-transitive for $n\geq2$ if $\mathcal{T}$ is 2-transitive. \\
\indent
Moreover, since $T=I_{\mathcal{K}}\otimes A$ commutes with $\mathcal{A}\otimes M_{n}(\mathbb{C})$, it means $A=\lambda I_{\mathbb{C}^{n}}$ because $A$ commutes $M_{n}(\mathbb{C})$ for some $\lambda\in\mathbb{C}$. So, $T=\lambda I_{\mathcal{K}\otimes\mathbb{C}^{n}}$, that is, $\mathcal{T}$ is 2-transitive by Arveson's Lemma 3.1. Then $\mathcal{T}$ is strong dense in $B(\mathcal{K}\otimes M_{n}(\mathbb{C}))$, that is $\mathcal{A}\otimes M_{n}(\mathbb{C})$ is catalytic.
\end{proof}
\indent Note that we can use the same proof to show that if $\mathcal{A}$ is a catalytic algebra on a Hilbert space $\mathcal{K}$, then $\mathcal{A}\otimes B(\mathcal{K}_{1})$ is also a catalytic algebra for any separable Hilbert space $\mathcal{K}_{1}$. Actually, the proof shows that $A\otimes K(\mathcal{K}_{1})$ is catalytic, where $K(\mathcal{K}_{1})$ denotes the algebra of compact operators. A more interesting question is whether $A\otimes I_{\mathcal{K}_{1}}$ is catalytic.\\
\indent
We have the following theorem by combining Theorem 4.4 and Proposition 4.5.
\begin{theorem}
For the flat unitary bundle $E$ over $R$, let $\mathcal{A}$ be a transitive subalgebra of $B(H^{2}_{E}(R))$ containing the algebra $\mathcal{J}_{E}(R)$. Then $\mathcal{A}$ is strongly dense in $B(H_{E}^{2}(R))$ or is a catalytic algebra.
\end{theorem}
\begin{proof}
Only one thing needs to be observed after noting that the bundles $R\times\mathbb{C}^{n}$ and $E$ extend to a trivial and a flat unitary bundle over the closure of $R$ and these extensions are similar [3]. If $\Phi$ is a bundle map from $clos(R)\times\mathbb{C}^{n}$ to the extension of $E$, establishing the similarity of $T_{\mathcal{H}}$ and $T_{E}$, then $\Phi$ induces a module isomorphism between $H^{2}_{\mathbb{C}^{n}}(R)$ and $H^{2}_{E}(R)$ conjugating $\mathcal{J}_{E}(R)$ and $\mathcal{J}_{\mathcal{R}\times\mathbb{C}^{n}}(R)$. Thus the similarity not only takes $T_{\mathcal{H}}$ to $T_{E}$, but also $\mathcal{J}(R)\otimes M_{n}(\mathbb{C})$ to $\mathcal{J}_{E}(R)$, which completes the proof.
\end{proof}
\indent
Note that for a separable Hilbert space $\mathcal{H}$, the proof goes through establishing that $\mathcal{J}(T_{E})$ and $M_{0}(E)$ are catalytic, where $\Phi\in M_{0}(E)$ if $\Phi(z)\in K(E_{z})$ for $z\in R$.

\bibliographystyle{amsplain}

\end{document}